\documentclass[12pt]{article}

\usepackage{amsmath,amsfonts,amssymb}

\newtheorem{thm}{Theorem}
\newtheorem{lem}[thm]{Lemma}
\newtheorem{prop}[thm]{Proposition}
\newcommand{\QED}{\mbox{$\square$}\\[-.6em]}

\newcommand{\sh}[1]{{#1}}
\newcommand{\RR}{\mathbb{R}}
\newcommand{\CC}{\mathbb{C}}
\newcommand{\ZZ}{\mathbb{Z}}
\newcommand{\NN}{\mathbb{N}}
\newcommand{\QQ}{\mathbb{Q}}

\newcommand{\al}{\alpha}
\newcommand{\lam}{\lambda}
\newcommand{\Lam}{\Lambda}
\newcommand{\del}{\delta}
\newcommand{\Del}{\Delta}

\newcommand{\sig}{\sigma}
\newcommand{\Sig}{\Sigma}
\newcommand{\om}{\varpi}
\newcommand{\ep}{\epsilon}

\newcommand{\g}{\mathfrak{g}}
\newcommand{\h}{\mathfrak{h}}

\newcommand{\hs}{\h^{\!*}}
\newcommand{\ghat}{\widehat{\g}}
\newcommand{\hhat}{\widehat{\h}}
\newcommand{\nhat}{\widehat{\mathfrak n}}
\renewcommand{\v}{^{\!\vee}}

\newcommand{\What}{\mathop{\widehat{W}}\nolimits}
\newcommand{\Wtil}{\widetilde{W}}

\newcommand{\B}{{\mathcal{B}}}
\newcommand{\Bhat}{{\hat{\B}}}

\newcommand{\wt}{\mathop{\rm wt}\nolimits}

\newcommand{\Stab}{\mathop{\rm Stab}\nolimits}

\newcommand{\st}{\star}
\newcommand{\sigb}{\bar\sigma }

\newcommand{\ot}{{\vdash}}

\begin{document}

\title{Littelmann Paths\\
 for the Basic Representation\\
 of an Affine Lie Algebra}
\author{Peter Magyar}
\date{August 13, 2003}

\maketitle

\abstract{Let $\g$ be a complex simple Lie algebra, and $\ghat$
the corresponding untwisted affine Lie algebra.  Let $\hat
V(\Lam_0)$ be the basic irreducible level-one representation of
$\ghat$, and $\hat V_{\lam\v}(\Lam_0)$ the Demazure module
corresponding to the translation $-\lam\v$ in the affine Weyl
group.  Suppose $\lam\v$ is a sum of minuscule coweights of $\g$
(which exist if $\g$ is of classical type or $E_6, E_7$).

We give a new model for the crystal graphs of $\hat V(\Lam_0)$ and
$\hat V_{\lam\v}(\Lam_0)$ which combines Littelmann's path model
and the Kyoto path model.  As a corollary, we prove that
$V_{\lam\v}(\Lam_0)$ is isomorphic as a $\g$-module to a tensor
product of fundamental representations of $\g$. }

\section{Main Results}

\subsection{Product Theorems}

Let $\g$ be a complex simple Lie algebra and $\ghat$ the
corresponding untwisted affine Kac-Moody algebra.  The basic
representation $\hat V(\Lam_0)$ is the simplest and most important
$\ghat$-module (see Sec.~2.1 for definitions, as well as \cite[Ch.~14]{Kac},\cite[Ch.~10]{Pressley-Segal}).  \
One of its remarkable properties is the Tensor Product Phenomenon.
In many cases, the Demazure modules $\hat V_z(\Lam_0)\subset \hat
V(\Lam_0)$ are representations of the finite-dimensional algebra
$\g$, and they factor into a tensor product of many small
$\g$-modules. Hence the full $\hat V(\Lam_0)$ could be constructed
by extending the $\g$-structure on the semi-infinite tensor power
$V\otimes V\otimes\cdots$ of a small $\g$-module $V$.

The Kyoto school of Jimbo, Kashiwara et
al., has
established this phenomenon for $\g$ of classical type (and for a
large class of $\ghat$-modules $\hat V(\Lam)$) via the theory of
perfect crystals \cite{KKMMNN},\cite{KMOTU1},\cite{KMOTU2},\cite{Hong-Kang},\cite{Kashiwara2} a development of the earlier theory of
semi-infinite paths \cite{DJKMO}. Pappas and Rapoport \cite{Pappas-Rapoport} have
given a geometric version of the phenomenon for type $A$:  they
construct a flat deformation of Schubert varieties of the affine
Grassmannian into a product of finite Grassmannians.

In this paper, we extend the Tensor Product Phenomenon for $\hat
V(\Lam_0)$ to the non-classical types $E_6$ and $E_7$ by a uniform
method which applies whenever $\g$ possesses a minuscule
representation, or more precisely a minuscule coweight. We shall
rely on a key property of such coweights which may be taken as the definition.  Let
$\hat X$ be the extended Dynkin diagram (the diagram of $\ghat$).
A coweight $\om\v$ of $\g$ is {\it minuscule} if and only if it is
a fundamental coweight $\om\v=\om_i\v$ and there exists an
automorphism $\sig$ of $\hat X$ taking the node $i$ to the
distinguished node 0. Such automorphisms exist in types
$A,B,C,D,E_6,E_7$.

We let $V(\lam)$ denote the irreducible $\g$-module with highest
weight $\lam$, and $V(\lam)^*$ its dual module.  Our main
representation-theoretic result is:
\begin{thm}
Let  $\lam\v$ be an element of the coroot lattice of $\g$ which is
a sum:
$$
\lam\v=\lam_1\v+\cdots+\lam_m\v,
$$
where $\lam_1\v,\cdots,\lam_m\v$ are minuscule fundamental
coweights (not necessarily distinct), with corresponding
fundamental weights $\lam_1,\cdots,\lam_m$.  Let $\hat
V_{\lam\v}(\Lam_0)\subset \hat V(\Lam_0)$ be the Demazure module
corresponding to the anti-dominant translation $t_{-\lam\v}$ in
the affine Weyl group.

Then there is an isomorphism of $\g$-modules:
$$
\hat V_{\lam\v}(\Lam_0)\cong V(\lam_1)^*\otimes\cdots\otimes
V(\lam_m)^*\,.
$$
\end{thm}

Now fix a minuscule coweight $\om\v$ and its corresponding
fundamental weight $\om$.  Let $N$ be the smallest positive
integer such that $N\om\v$ lies in the coroot lattice of $\g$.
Then we have the following characterization of the basic
irreducible $\ghat$-module:
\begin{thm} The tensor power $V_N{:=}V(\om)^{\otimes N}$ possesses
non-zero $\g$-invariant vectors.  Fix such a vector $v_N$, and
define the $\g$-module $V^{\otimes\infty}$ as the direct limit of
the sequence:
$$
V_N\hookrightarrow V_N^{\otimes 2}\hookrightarrow V_N^{\otimes
3}\hookrightarrow\cdots
$$
where each inclusion is defined by: $v\mapsto v_N\otimes v$.

Then $\hat V(\Lam_0)$ is isomorphic as a $\g$-module to
$V^{\otimes\infty}$.
\end{thm}
It would be interesting to define the action of the full algebra $\ghat$ on $V^{\otimes\infty}$, and thus
give a uniform ``path construction'' of the basic representation
(cf.~\cite{DJKMO}):  that is, to define the raising and lowering operators $E_0, F_0$, as well as the charge operator $d$.  Combinatorial definitions of the charge produce for $\g$ of classical type produce generalizations of the Kostka-Foulkes polynomials (c.f. \cite{Okado}).

\subsection{Crystal Theorems}

Our basic tool to prove the above results is Littelmann's
combinatorial model
\cite{Littelmann1},\cite{Littelmann2},\cite{LLM2} for
representations of Kac-Moody algebras, a vast generalization of
Young tableaux.  Littelmann's paths and path operators give a
flexible construction of the crystal graphs associated to quantum
$\g$-modules by Kashiwara \cite{Kashiwara1}
and Lusztig \cite{Lusztig} (see also \cite{Joseph},\cite{Hong-Kang}).  Roughly speaking, we prove Theorem 1 by reducing
it to an identity of paths: we construct a path crystal for the
affine Demazure module which is at the same time a path crystal
for the tensor product. We carry this out in Sec.~2.

Theorem 2 follows as a corollary in Sec.~3.  To describe the
crystal graph of the semi-infinite tensor product, we pass to a
semi-infinite limit of Littelmann paths which we call {\it
skeins}.  We thus recover the Kyoto path model for classical $\g$,
and our results are equally valid for $E_6$, $E_7$.

To be more precise, we briefly sketch Littelmann's theory. Let us
define a {\it $\g$-crystal} as a {\it set} $\B$ with a {\it weight
function} $\wt:\B\to\oplus_{i=1}^r\ZZ\om_i$, and partially defined
{\it crystal operators}\ $e_1,\ldots,e_r$,$f_1,\ldots,f_r:\B\to
\B$ satisfying:
$$
\wt(f_i(b))=\wt(b)-\al_i\qquad\text{and}\qquad
e_i(b)\sh=b'\Longleftrightarrow f_i(b')\sh=b\,.
$$
Here $\om_1,\ldots,\om_r$ are the fundamental weights and
$\al_1,\ldots,\al_r$ are the roots of $\g$.  A {\it dominant}
element is a $b\sh\in B$ such that $e_i(b)$ is not defined for any
$i$. We say that a crystal $\B$ is a {\it model} for a $\g$-module
$V$ if the formal character of $\B$ is equal to the character of
$V$, and the dominant elements of $\B$ correspond to the
highest-weight vectors of $V$. That is:
$$
\textstyle\mathop{\rm char}(V)=\sum_{b\in B}e^{\wt(b)}
\qquad\text{and}\qquad V\cong\bigoplus_{b\ \text{dom}}
V(\wt(b))\,,
$$
where the second sum is over the {dominant} elements of $\B$.
Clearly, a $\g$-module $V$ is determined up to isomorphism
by any model $\B$.

We construct such $\g$-crystals $\B$ consisting of polygonal paths
in the vector space of weights, $\hs_\RR:=\oplus_{i=1}^r\RR\om_i$.
Specifically:

\begin{itemize}

\item The {\it elements} of $\B$ are continuous piecewise-linear
mappings $\pi:[0,1]\to\hs_\RR$, up to reparametrization, with
initial point $\pi(0)=0$. We use the notation $\pi=(v_1\st
v_2\st\cdots\,\st v_k)$, where $v_1,\ldots,v_k\in\hs_\RR$ are
vectors, to denote the polygonal path starting at 0 and moving
linearly to $v_1$, then to $v_1\sh+v_2$, etc.

\item The {\it weight} of a path is its endpoint:
$$
\wt(\pi) := \pi(1) = v_1\sh+\cdots\sh+v_k\,.
$$

\item The {\it crystal lowering operator} $f_i$ is defined as
follows (and there is a similar definition of the raising operator
$e_i$). Let $\st$ denote the natural associative operation of
concatenation of paths, and let any linear map
$w:\hs_\RR\to\hs_\RR$ act pointwise on paths: $w(\pi):=
(w(v_1)\st\cdots\st w(v_k))$. We will divide a path $\pi$ into
three well-defined sub-paths, $\pi=\pi_1\st\pi_2\st\pi_3$, and
reflect the middle piece by the simple reflection $s_i$:
$$
f_i\pi:=\pi_1\st s_i\pi_2\st\pi_3\,.
$$

\qquad The pieces $\pi_1,\pi_2,\pi_3$ are determined according to
the behavior of the $i$-height function
$h_i(t)=h_i^\pi(t):=\langle\pi(t),\al_i\v\rangle$, as the point
$\pi(t)$ moves along the path from $\pi(0)=0$ to
$\pi(1)=\wt(\pi)$. This function may attain its minimum value
$h_i(t)=M$ several times. If, after the {\it last} minimum point,
$h_i(t)$ never rises to the value $M\sh+1$, then $f_i\pi$ is {\it
undefined}. Otherwise, we define $\pi_2$ as the last sub-path of
$\pi$ on which $M\leq h_i(t)\leq M\sh+1$, and $\pi_1$, $\pi_3$ as
the remaining initial and final pieces of $\pi$.

\end{itemize}

A key advantage of the path model is that the crystal operators,
while complicated, are universally defined for all paths.  Hence a
path crystal is completely specified by giving its set of paths
$\B$.

Also, the dominant elements have a neat pictorial
characterization, as the paths $\pi$ which never leave the
fundamental Weyl chamber: that is, $\h_i^\pi(t)\geq 0$ for all
$t\in[0,1]$ and all $i=1,\ldots,r$.  For simplicity we restrict
ourselves to {\it integral} dominant paths, meaning that all the
steps are integral weights:
$v_1,\ldots,v_k\in\oplus_{i=1}^r\ZZ\om_i$.
(For arbitrary dominant paths, see \cite{Littelmann2}.)

Littelmann's Character Theorem \cite{Littelmann2} states
that if $\pi$ is any integral dominant path with weight $\lam$,
then the set of paths $\B(\pi)$ generated from $\pi$ by
$f_1,\ldots,f_r$ is a model for the irreducible $\g$-module
$V(\lam)$.  (This $\B(\pi)$ is also closed under
$e_1,\ldots,e_r$.) Note that we can choose {\it any} integral path
$\pi$ which stays within the Weyl chamber and ends at $\lam$, and
each such choice gives a different (but isomorphic) path crystal
modelling $V(\lam)$. In principle, any reasonable indexing set for
a basis of $V(\lam)$ should be in natural bijection with $\B(\pi)$
for some choice of $\pi$.  For example, classical Young tableaux
correspond to choosing the steps $v_j$ to be coordinate vectors in
$\hs_\RR\cong\RR^n$.

Furthermore, we have Littelmann's Product Theorem \cite{Littelmann2}: if $\pi_1,\ldots,\pi_m$ are dominant integral
paths of respective weight $\lam_1$,\ldots, $\lam_m$, then
$\B(\pi_1)\st\cdots\st\B(\pi_m)$, the set of all concatenations,
is a model for the tensor product $V(\lam_1)\otimes\cdots\otimes
V(\lam_m)$.

Everything we have said also holds for the affine algebra $\ghat$,
provided we replace the roots $\al_1,\ldots,\al_r$ of $\g$ by the
roots $\al_0,\al_1,\ldots,\al_r$ of $\ghat$ ; and the weights
$\om_1,\ldots,\om_r$ of $\g$ by the weights
$\Lam_0,\Lam_1,\ldots,\Lam_r$ of $\ghat$.  We also replace the
vector space $\hs_\RR$ by $\hhat^*_\RR:=
\oplus_{i=0}^r\RR\Lam_i\,\oplus\,\RR\del$, where $\del$ is the
non-divisible positive imaginary root of $\ghat$. (Indeed, the
theory works uniformly for all symmetrizable Kac-Moody algebras.)
We denote path crystals for $\g$ and $\ghat$ by $\B$ and $\Bhat$
respectively.

We can also model the affine Demazure module $\hat V_z(\Lam)$,
where $z\in\What$, the Weyl group of $\ghat$.  Indeed, if
$z=s_{i_1}\cdots s_{i_m}$ is a reduced decomposition and $\pi$ is
an integral dominant path of weight $\Lam$, we define the {\it
Demazure path crystal}:
$$
\Bhat_z(\pi):=\{f_{i_1}^{k_1}\cdots f_{i_m}^{k_m}\pi\,\mid\,
k_1,\ldots,k_m\geq 0\}\,.
$$
Then the formal character of $\Bhat_z(\pi)$ is equal to the
character of $\hat V_z(\Lam)$, and $\pi$ is the unique dominant path \cite{Littelmann1}.  Now suppose $z=t_{-\lam\v}$, an
anti-dominant translation in $\What$, so that $\hat
V_{\lam\v}(\Lam):=\hat V_z(\Lam)$ is a $\g$-submodule of $\hat
V(\Lam)$; and consider $\Bhat_{\lam\v}(\pi):=\Bhat_z(\pi)$ as a
$\g$-crystal by forgetting the action of $f_0,e_0$ and projecting
the affine weight function to $\hs_\RR$.  Then Littelmann's
Restriction Theorem \cite{Littelmann2} implies that the
$\g$-crystal $\Bhat_{\lam\v}(\pi)$ is a model for the $\g$-module
$\hat V_{\lam\v}(\Lam)$.

Now we are ready to state our main combinatorial results.  For
$\lam$ a dominant weight, define its dual weight $\lam^{\!*}$ by
the dual $\g$-module: $V(\lam^{\!*})=V(\lam)^*$.
\begin{thm}
Let $\lam\v$ be as in Theorem 1. Then there exists an integral
$\ghat$-dominant path $\pi$ with weight $\Lam_0$ which generates
the $\ghat$ Demazure path crystal:
$$
\Bhat_{\lam\v}(\pi)= \Lam_0\st\B(\lam_1^*)\st\cdots\st\B(\lam_m^*)
\!\!\!\mod\RR\del\,.
$$
This is to be understood as an equality of sets of paths in
$\hhat^*_\RR\!\mod\RR\del$.
\end{thm}
Theorem 1 follows immediately from this.  Indeed,
$s_i\Lam_0=\Lam_0$ for $i=1,\ldots,r$, so
$f_i(\Lam_0\st\pi')=\Lam_0\st f_i(\pi')$ for any path $\pi'$. Thus
the crystal on the right-hand side of the above equation is
isomorphic to $\B(\lam_1^*)\st\cdots\st\B(\lam_m^*)$, which models
$V(\lam_1)^*\otimes\cdots\otimes V(\lam_r)^*$.
See \cite{Grabiner-Magyar} for methods of enumerating the paths in this crystal (and hence computing the dimension of the corresponding representation).  

Next we give the following crystal version of Theorem 2:
\begin{thm}
Let $\om\v$,\, $N$ be as in Theorem 2. Define the $N$-fold
concatenation $\B_N=\B(\om^*)\st\cdots\st\B(\om^*)$.  Then
$\Lam_0\st\B_N$ contains a unique $\ghat$-dominant path
$\Lam_0\st\pi_N$.

Define $\Bhat_\infty$ as the direct limit of the sequence:
$$
\Lam_0\st\B_N\stackrel{}\hookrightarrow\Lam_0\st
\B_N\st\B_N\stackrel{}\hookrightarrow\Lam_0\st
\B_N\st\B_N\st\B_N \stackrel{}\hookrightarrow\cdots\,
$$
where the inclusions are given by $\Lam_0\sh\st\pi\mapsto \Lam_0\sh\st\pi_N\sh\st\pi$.
Then $\Bhat_\infty$ has a natural $\ghat$-crystal structure which is isomorphic to the
$\ghat$-crystal of $\hat V(\Lam_0)$.
\end{thm}
In Sec.~3.2, we will use crystals of semi-infinite paths (``skeins") to give a natural meaning to the direct limit above.

\subsection{Example: $E_6$}
Referring to Bourbaki \cite{Bourbaki}, we write the extended
Dynkin diagram $\hat X=\hat E_6$:

$$\begin{array}{c@{\!}c@{\!}c@{\!}c@{\!}c@{\!}c@{\!}c@{\!}c@{\!}c}
&&&&\circ&\mbox{\scriptsize 0}&&&\\[-.55em]
&&&&|&&&&\\[-.55em]
&&&&\bullet&\mbox{\scriptsize 2}&&&\\[-.55em]
&&&&|&&&&\\[-.55em]
\bullet&\mbox{---}&\bullet&\mbox{---}&\bullet&\mbox{---}&
\bullet&\mbox{---}&\bullet\\[-.5em]
\mbox{\scriptsize 1}&&\mbox{\scriptsize  3}&&\mbox{\scriptsize 4}&
&\mbox{\scriptsize  5}&&\mbox{\scriptsize  6}
\end{array}$$
\\
The simple roots are defined inside $\RR^6$ with standard basis
$\ep_1,\ldots,\ep_6$. (Our $\ep_6$ is $\frac{1}{\sqrt
3}(-\ep_6\sh-\ep_7\sh+\ep_8)$ in Bourbaki's notation.)  They are:
$$\begin{array}{c}
\al_1 =\tfrac12(\ep_1\sh+\ep_2\sh+\ep_3\sh+\ep_4+\ep_5)\sh+\tfrac{\sqrt 3}{2}\ep_6,\quad\al_2=\ep_1\sh+\ep_2,\\[.2em]
\al_3=\ep_{2}\sh-\ep_{1},\quad\al_4=\ep_{3}\sh-\ep_{2}
,\quad\al_5=\ep_{4}\sh-\ep_{3},\quad\al_6=\ep_{5}\sh-\ep_{4}\,.
\end{array}$$
Since $E_6$ is simply laced, the coroots and coweights may be
identified with the roots and weights, with the natural pairing
given by the standard dot product on $\RR^6$.

We focus on the minuscule coweight $\om_1\v$ corresponding to the
diagram automorphism $\sig$ with $\sig(1)=0$ and $\sig(0)=6$. In
this case, the corresponding fundamental representation $V(\om_1)$
is also minuscule, meaning that all of its weights are extremal
weights $\lam\in W(E_6)\sh\cdot\om_1$. The roots
$\al_2,\cdots,\al_6$ generate the root sub-system $D_5\subset
E_6$, and the reflection subgroup $W(D_5)=\Stab_{W(E_6)}(\om_1)$
acts by permuting $\ep_1,\ldots,\ep_5$ (the subgroup $W(A_4)=S_5$) and
by changing an even number of signs $\pm\ep_1,\ldots,\pm\ep_5$.  We
have $\dim V(\om_1)={|W(E_6)/W(D_5)|}=27$. The weights are:
$$\begin{array}{c}
\om_1\sh=\tfrac{2\sqrt3}3\ep_6,\\[.3em]
S_5\sh\cdot\tfrac12(-\ep_1\sh+\ep_2\sh+\ep_3\sh+\ep_4\sh+\ep_5)
+\tfrac{\sqrt3}6\ep_6,\\[.3em]
S_5\sh\cdot\tfrac12(-\ep_1\sh-\ep_2\sh-\ep_3\sh+\ep_4\sh+\ep_5)
+\tfrac{\sqrt3}6\ep_6,\\[.3em]
-\tfrac12(\ep_1\sh+\ep_2\sh+\ep_3\sh+\ep_4\sh+\ep_5)
+\tfrac{\sqrt3}6\ep_6,\\[.3em]
\pm S_5\sh\cdot\ep_1 -\tfrac{\sqrt3}3\ep_6,
\end{array}$$
The lowest weight is $-\om_6=-\ep_5\sh-\tfrac{\sqrt3}3\ep_6$, so
that $V(\om_1)^*=V(\om_6)$ and $\om_1^*=\om_6$.

The simplest path crystal for $V(\om_1^*)$ is the set of 27
straight-line paths from 0 to the negatives of the above extremal
weights:
$$
\B(\om_1^*)=\{\,(v)\,\mid\, v\in -W(E_6)\sh\cdot\om_1\,\}
$$
We have $3\om_1\v\in\oplus_{i=1}^6\RR\al_i\v$ the coroot lattice,
so that $N=3$ in Theorem 2, and this $N$ is also the order of the
automorphism $\sig$.  The path crystal
$\B_3:=\B(\om_1^*)\st\B(\om_1^*)\st\B(\om_1^*)$, the set of all
3-step walks with steps chosen from the 27 weights of
$V(\om_1^*)$, is a model for $V(\om_1^*)^{\otimes 3}$. In this
case there is a unique $\g$-dominant path of weight 0,
$$
\pi_3:=(\om_6)\st(\om_1\sh-\om_6)\st(-\om_1)\,,
$$
which corresponds to the one-dimensional space of $\g$-invariant
vectors in $V(\om_i^*)^{\otimes 3}$.

Now Theorem 3 states that the affine Demazure module $\hat
V_{3m\om_1\v}(\Lam_0)$ is modelled by the $\ghat$-path crystal:
$$
\B_{3m}=\{(\Lam_0\st v_1\st\cdots\st v_{3m})\mid
v_j\in-W(E_6)\sh\cdot\om_1\}\,,
$$
the set of all $3m$-step walks in $\Lam_0\oplus\RR^6$ starting
at $\Lam_0$, with steps chosen from the 27 weights of
$V(\om_1^*)$.  This path crystal is generated from its unique
$\ghat$-dominant path $\Lam_0\sh\st\pi_3\sh\st\cdots\sh\st\pi_3$. As a
corollary of Theorem 4, we can realize the $\ghat$-crystal of the
basic $\ghat$-module $\hat V(\Lam_0)$ as the set of all infinite
walks (or ``skeins'') of the form:
$$
\pi=\Lam_0\st\underbrace{\pi_3\st\cdots\st\pi_3}_{\text{infinite}}\st\,
v_1 \st\cdots\st v_{3m}\,,
$$
with $m>0$ and $v_j\in-W(E_6)\sh\cdot\om_1$.  The endpoint of such
a skein is $\wt(\pi):=\Lam_0\sh+v_1\sh+\cdots\sh+v_m$. The crystal
operators $f_i$ are defined just as for finite paths.  Acting near the
end of the skein, they unwind the coils $\pi_3$ one at a time,
right-to-left. See Sec.~3 for details.

\section{Demazure Crystals}

\subsection{Notations}

We will work with a complex simple Lie algebra $\g$ of rank $r$, a
Cartan subalgebra $\h\subset\g$, the set of roots
$\Del\subset\hs$\!,  and  the set of coroots $\Del\v\subset\h$. We
write the highest root of $\Del$ as
$\theta=a_1\al_1+\cdots+a_r\al_r$, and its coroot as
$\theta\,\v=a_1\v\al_1\v+\cdots+a_r\v\al_r\v$. Warning: If $\g$ is
not simply laced, $\theta\,\v$ is not the highest root of the dual
root system $\Del\v$.

The Weyl group $W$ of $\g$ is generated by reflections
$s_1,\ldots,s_r$ defined by
$s_i(\lam)=\lam-\langle\lam,\al_i\v\rangle\al_i$ for
$\lam\in\hs_\RR:=\oplus_{i=1}^r\RR\om_i$.  We have the fundamental
Weyl chamber $C=\{\lam\in\hs_\RR\mid
\langle\lam,\al_i\v\rangle\sh\geq0,\ \ i\sh=1,\ldots,r\}$. The
Weyl group also acts naturally on $\h_\RR$.  If we choose a
$W$-invariant bilinear form $(\,\mathbf\cdot\,|\,\mathbf\cdot\,)$
on $\h_\RR$, we have the isomorphism $\nu:\h_\RR\to\hs_\RR$
defined by $\langle\nu(h),h'\rangle=(h|h')$ for $h,h'\in\h_\RR$.
We normalize so that $\nu(\theta\,\v)=\theta$ and
$\nu(\om_i\v)=\frac{a_i}{a_i\v}\om_i$.

Now let $\ghat=\g\otimes\CC[t,t^{-1}]\,\oplus\,\CC K\,\oplus\,\CC
d$ be the untwisted affine Lie algebra of $\g$, where $K$ is a
central element and $d=t\frac{d}{dt}$ is a derivation.  (Cf.~Kac
\cite[Ch.~6 and 7]{Kac}.)  Then $\ghat$ has Cartan subalgebra
$\hhat=\h\oplus\CC K\oplus\CC d$, with dual
$\hhat^*=\hs\oplus\CC\Lam_0\oplus\CC\del$, where
$\langle\Lam_0,\h\rangle=\langle\del,\h\rangle=0$ and
$\langle\Lam_0,K\rangle=\langle\del,d\rangle=1$.

The simple roots of $\ghat$ are $\al_1,\ldots,\al_r$ and
$\al_0=\del\sh-\theta$; the simple coroots are
$\al_1\v,\ldots,\al_r\v$ and $\al_0\v=K\sh-\theta\,\v$. The
fundamental weights are $\Lam_0$ and $\Lam_i=\om_i+a_i\v\Lam_0$
for $i=1,\ldots,r$. The affine Weyl group $\What$ is generated by
the reflections $s_0,s_1,\ldots,s_r$ acting on $\hhat^*_\RR$. The
fundamental Weyl chamber of $\ghat$ is the cone
$\hat{C}=\{\mbox{$\Lam\in\hhat^*_\RR$}\mid \langle
\Lam,\al_i\v\rangle\geq 0,\ i\sh=0,\ldots,r\}$ with extremal rays
$\Lam_0,\ldots,\Lam_r$.

For $\Lam\in\oplus_{i=0}^r\NN\Lam_i$ and $w\in \What$, we have the
irreducible highest-weight $\ghat$-module $\hat V(\Lam)$.  We will
also consider the Demazure module $\hat V_z(\Lam):=U(\nhat_+)
\sh\cdot v_{z\Lam}$, where $\nhat_+$ is the algebra spanned by the
positive weight-spaces of $\ghat$,\ $z\in \What$ is a Weyl group
element, and $v_{z\Lam}$ is a non-zero vector of extremal weight
$z\Lam$ in $\hat V(\Lam)$.

For a weight $\Lam=\lam+\ell\Lam_0+m\del$, we will ignore the
charge $m=\langle\Lam,d\rangle$ and work modulo $\RR\del$, even
when we do not indicate this explicitly. Since
$\langle\del,\al_i\v\rangle=0$ for $i=0,\ldots,r$, the charge has no effect on the
path operators $f_i,e_i$.

Consider the lattice $M=\nu(\oplus_{i=1}^r\ZZ\al_i\v)$ in
$\hs_\RR$.  For any $\mu\in M$, there is an element
$t_\mu\in\What$ which acts on the weights of level $\ell$ as
translation by $\ell\mu$: that is,
$$
t_\mu(\lam\sh+\ell\Lam_0)=\lam\sh+\ell\mu\sh+\ell\Lam_0\!\pmod{\RR\del}.
$$
Furthermore, the affine Weyl group is a semi-direct product of the
finite Weyl group with the lattice of translations: $\What =
W\ltimes t_M$.

Consider the anti-dominant translation $z=t_{-\lam}$ corresponding
to a dominant weight $\lam=\nu(\lam\v)\in C\cap M$. We denote the
resulting Demazure module as $\hat V_{\lam\v}(\Lam_0):=\hat
V_{z}(\Lam_0)$. Then the $\nhat_+$-module $\hat
V_{\lam\v}(\Lam_0)$ is also a $\g$-submodule of $\hat V(\Lam_0)$:
$$
\g\sh\cdot\hat V_{\lam\v}(\Lam_0)\subset\hat V_{\lam\v}(\Lam_0)\,
,
$$
and these are the only $z\in\What$ for which $\hat V_z(\Lam_0)$ is
a $\g$-module.

\subsection{Minuscule weights and coweights}

We collect needed facts concerning minuscule weights in root
systems. The statements below are well-known and easily verified
from tables \cite{Bourbaki},\cite{Kac}, although direct proofs are
also not difficult (cf. \cite{Macdonald}).

We say a non-zero coweight $\om\v\in\hs_\RR$ is {\it minuscule for
$\Del$} if $\langle\al,\om\v\rangle=0$ or 1 for all positive roots
$\al\in\Del_+$.  Equivalently, $\om\v=\om_i\v$ for some
$i=1,\ldots,r$ with $a_i=\langle\theta,\om_i\v\rangle=1$. This
implies that $a_i\v=1$ as well, so that $\nu(\om_i\v)=\om_i$. The
classification of the minuscule $\om_i\v$ is most concisely
described by listing the pairs $(X,\,X\sh\setminus\!\{i\})$, where
$X$ is the Dynkin diagram of $\g$.  We have $\om_i\v$ minuscule
when:
$$\begin{array}{rcl}
(X,\,X\sh\setminus\!\{i\})&\cong&(A_r,A_{r-k}\sh\times A_{k-1}),\
\ k\sh=1,\ldots,r\\&& (B_r,B_{r-1}),\ \ (C_r,A_{r-1}),\\&&
(D_r,D_{r-1}),\ \ (D_r,A_r),\\&& (E_6,D_5),\ \ (E_7,E_6)\,.
\end{array}
$$
There are no minuscule $\om_i\v$ for $X=E_8, F_4$, or $G_2$.

Now define the extended Weyl group $\Wtil$ as a group of linear
mappings on $\hhat^*_\RR$:  namely, $\Wtil:=W\ltimes t_L$, where
$L=\nu(\oplus_{i=1}^r\ZZ\om_i\v)$.
Let
$$\Sig:=\{\sig\in\Wtil\mid \sig(\hat{C})=\hat{C}\}\,,$$ the
symmetries in $\Wtil$ of the fundamental chamber of $\hhat^*_\RR$.
The set $\Sig$ is a system of coset representatives for
$\Wtil/\What$, so that $\Wtil=\Sig\ltimes\What$. We can extend the
Bruhat length function to $\Wtil$ as: $l(\sig w)=l(w\sig):=l(w)$
for $\sig\in\Sig$, $w\in\What$.
Each element $\sig\in\Sig$ defines an automorphism of the Dynkin
diagram of $\ghat$ which we also write as $\sig$. For
$j=0,\ldots,r$, we have:
$$
\sig(\Lam_j)=\Lam_{\sig(j)}\quad\text{and}\quad
\sig(\al_j)=\al_{\sig(j)}\,.
$$

There is a natural correspondence between elements of $\Sig$ and
minuscule coweights.  Each $\sig\in\Sig$ can be written uniquely
as:
$$
\sig=\sigb\,t_{-\nu(\om_i\v)}=\sigb\,t_{-\om_i}\,,
$$
for $\sigb\in W$ and $\om_i\v$ a minuscule coweight. In fact,
$\sigb=w_0w_i$, where $w_0$ is the longest element of $W$ and
$w_i$ is the longest element of the parabolic subgroup
$W_i:=\Stab_W(\om_i)$.  We have
$\sigb(\al_j)=\al_{\sig(j)}$ for $j\neq i$, and
$\sigb(\al_i)=-\theta$.

We have $\sig(\Lam_i)=\sigb \mbox{$t_{-\om_i}(\om_i\sh+\Lam_0)$}$
$=\Lam_0$, so that
$$\sig(i)=0\,.$$
Also, $\sigb(\om_i)=w_0w_i(\om_i)=w_0(\om_i)=-\om_i^*$, and
$\om_{\sig(0)}=\om_i^*$.

The number $N$ appearing in Theorems 2 and 4 is the order of $\sig$ in the group $\Sig$, and is also the order of $\om\v$ in $\oplus_{i=1}^r\RR\om_i\v\,/\,\oplus_{i=1}^r\RR\al_i\v$.

The definition of a minuscule weight $\om$ for $\Del\v$ is dual to
the above: $\langle\om,\al\v\rangle=0$ or 1 for all positive
coroots $\al\v\in\Del\v_+$; or equivalently $\om=\om_i$ and
$\langle\om,\theta^*\rangle=1$, where $\theta^*$ is highest in the
root system $\Del\v$.  The fundamental representation $V(\om)$ corresponding to a minuscule $\om$ has a basis consisting of extremal weight vectors $v_{w(\om)}$ for $w\in W$.  Note, however, that $\om$ need not be minuscule even when the corresponding coweight $\om\v$ is minuscule.
In fact, we have:
\begin{lem} Let $\om_i\v$ be a minuscule coweight for $\Del$,
and $\om_i$ the corresponding weight. 
\\[.3em]
(i) If $\Del$ is simply laced (i.e., all root vectors have the same length), then $\om_i$ is
minuscule for $\Del\v$.
\\[.3em]
(ii) If $\Del$ is not simply laced, then $\om_i$ is minuscule for
$\Del\v_s$, the simply-laced root system of short vectors in
$\Del\v$.  Furthermore, $\al_i\v\in\Del\v_s$.
\end{lem}
{\it Proof.}  Part (i) follows from
$\nu(\om_i\v)=(a_i/a_i\v)\,\om_i=\om_i$. Part (ii) is immediately
verified for the relevant types $B_r$ and $C_r$. \QED

Consider the {\it parabolic Bruhat order} on the $W$-orbit
$W\sh\cdot\om_i$.  That is, the partial order generated by the
relations: $\tau_1< \tau_2$ if $\tau_1=\tau_2-d\al$ for some
positive root $\al$ and some $d>0$. 
The following result says that
if $\om_i\v$ is minuscule, then this ``strong'' order is identical
to the ``weak'' order:

\begin{lem}[Stembridge \cite{Stembridge}]
Suppose the coweight $\om_i\v$ is minuscule. Then the Bruhat order
on $W\sh\cdot\om_i$ has covering relations:
$$\tau_1\lessdot \tau_2\ \
\text{whenever}\ \ \tau_1=\tau_2-d\al_j $$ for some \emph{simple}
root $\al_j$ of $W$ and some $d>0$.
\end{lem}
{\it Proof.}  This follows from Stembridge's formulation by noting
that $W\sh\cdot\om_i\cong W/W_i\cong W\sh\cdot\om_i\v$. \QED


\subsection{Lakshmibai-Seshadri Paths}

We examine in detail the path crystals of $V(\om_i)$ where
$\om_i\v$ is minuscule.  For any dominant weight $\lam$, the most
canonical choice of dominant path is the straight-line path from 0
to $\lam$, denoted $\pi=(\lam)$.  The corresponding path crystal
$\B(\lam)$ can be described non-recursively by the {\it
Lakshmibai-Seshadri (LS) chains} \cite{Littelmann1}. These are
saturated chains in the parabolic Bruhat order on $W\sh\cdot\lam$,
weighted with certain rational numbers:
$$
(\tau_1\gtrdot\cdots\gtrdot\tau_m;\ 0\sh<a_1\sh\leq\cdots\sh\leq
a_{m-1}\sh<1)
$$
with $\tau_j\sh\in W\sh\cdot\lam$,\ \ $a_j\sh\in\QQ$, and
$m\sh\geq1$. We require that if $\tau_{j+1}=\tau_j-d_j\al$ for
$\al\sh\in\Del_+$,\ \ $d_j\sh\in\NN$, then $a_j=n_j/d_j$ for some
$n_j\sh\in\NN$.  An LS chain corresponds to the {\it LS path}
defined as:
$$
\pi=(a_1\tau_1\st(a_2\sh-a_1)\tau_2\st(a_3\sh-a_2)\tau_3\st\cdots
\st(1\sh-a_{m-1})\tau_m)\,.
$$
Notice that if $a_{j+1}=a_j$ then we may omit the step $0\tau_j$;
in this case there may be more than one LS chain producing the
same LS path.  However, $B(\lam)$ is equal to the set of all
distinct LS paths \cite{Littelmann1}.


\begin{prop}
Let $\om_i\v, \om_l\v$ be two minuscule coweights (possibly
identical), and let $\sigb_l$ be the linear mapping of $\hs_\RR$
corresponding to $\om_l\v$.  For each $\pi\in\B(\om_i)$, we have
$\sigb_l\pi\in\B(\om_i)$. That is, the linear mapping $\sigb_l$
permutes the paths in $\B(\om_i)$.
\end{prop}
{\it Proof.}  Consider an LS chain
$(\tau_1\sh\gtrdot\cdots\sh\gtrdot\tau_m;\
0\sh<a_1\sh\leq\cdots\sh\leq a_{m-1}\sh<1)$ corresponding to a
path $\pi\in\B(\om_i)$. If $\g$ is simply laced, then $\om_i$ is
minuscule by Lemma 5(i), and the denominator of $a_j$ is:
$$
d_j=\langle\tau_{j},\al\v\rangle=\langle
w\om_i,\al\v\rangle=\langle \om_i,w\al\v\rangle=0\ \text{or}\ 1\,.
$$
Since $0\sh<a_j\sh=\frac{n_j}{d_j}\sh<1$, this means that $m=1$
and $\pi=(\tau_1)=(w\om_i)$ for $w\in W$, a straight-line path of
extremal weight.  Since $\sigb_l$ is an automorphism of the root
system $\Del$, it permutes the elements of a $W$-orbit, and hence
$\sigb_l\pi$ is another straight-line LS path.

If $\g$ is not simply laced, the paths of $\B(\om_i)$ are more
complicated.  By Lemma 6, we may assume that
$\tau_{j+1}=\tau_j-d_j\al_{k(j)}$, where $k(j)\in\{1,\ldots,r\}$
and $\al_{k(j)}$ is a {\it simple} root.  The turned path is:
$$
\sigb_l\pi=
(a_1\sigb_l\tau_1\st(a_2\sh-a_1)\sigb_l\tau_2\st\cdots)\,.
$$
Suppose $k(j)=l$ for some $j$. By Lemma 5(ii),\ $\al_l\v$ is a
short root of $\Del\v$, and so is $w\al_l\v$, so we have:
$$
d_j=\langle \tau_j,\al_l\v\rangle =\langle w\om_i\v,\al_l\v\rangle
=\langle \om_i,w\al_l\v\rangle=0\ \text{or}\ 1
$$
by the same Lemma. This again means that $m=1$ and $\sigb_l\pi$ is
a straight-line LS path.

Finally, suppose $k(j)\sh\neq l$ for all $j=1,\ldots,m$. Then:
$$
\sigb_l\tau_{j+1}=\sigb_l\tau_{j}-d_j\sigb_l\al_{k(j)}
=\sigb_l\tau_{j}-d_j\al_{p}\,,
$$
where $p:=\sig_l(k(j))\in\{1,\ldots,r\}$.  Hence
$\sigb_l\tau_{j}\gtrdot s_{p}\sigb_l\tau_j=\sigb_l\tau_{j+1}$, and
$\sigb_l\pi$ is an LS path. \QED
\\[.3em]
Although we do not need it here, we note that for any minuscule
$\om_i\v$, the LS-paths of $B(\om_i)$ have at most two linear
pieces (cf.~\cite{Lakshmibai-Seshadri}).

\subsection{Twisted Demazure operators}

For a Weyl group element with reduced decomposition
$z=s_{i_1}\cdots s_{i_m}\in\What$ and any path $\pi$ (not
necessarily dominant), we define the {\it crystal Demazure
operator}:
$$
\Bhat_z(\pi):=\{f_{i_1}^{k_1}\cdots
f_{i_m}^{k_m}\pi_\lam\,\mid\, k_1,\ldots,k_m\geq 0\}\,.
$$
%
%
%
%
We can extend $\Bhat_z$ to an operator taking any set of paths
$\Pi$ to a larger set of paths: $\Bhat_z(\Pi):=\bigcup_{\pi\in\Pi}
\Bhat_z(\pi)$.   We have $
\Bhat_y(\Bhat_z(\Pi))=\Bhat_{yz}(\Pi)$ whenever $l(yz)=l(y)+l(z).$\
Similarly, we let $\B(\Pi)$ be the set of all paths generated from
$\Pi$ by $f_1,\ldots,f_r$, $e_1,\ldots,e_r$.

It will be convenient to define a Demazure module $\hat V_z(\Lam)$
for any $z\in\Wtil=\Sig\ltimes\What$.  Now, $\sig\in\Sig$ induces
an automorphism of $\ghat$, so for a module $\hat V$ we have the twisted
module $\sig \hat V$ defined by the action $g{\odot} v:=\sig^{-1}(g)v$
for $g\in\ghat$, $v\in \hat V$. That is, $\sig\hat V(\Lam)\cong\hat
V(\sig\Lam)$, and in particular $\sig\hat V(\Lam_i)\cong\hat
V(\Lam_{\sig(i)})$. Now, for $z=\sig y$ with $y\in\What$ define: $
\hat V_{z}(\Lam):=\sig(\hat V_y(\Lam))\subset \sig\hat V(\Lam)\,, $ a
twist of an ordinary Demazure module. Thus, $\hat V_{\sig
y}(\Lam)\cong \hat V_{\sig y\sig^{-1}}(\sig\Lam)$ and $$\hat
V_{y\sig}(\Lam)\cong \hat V_y(\sig\Lam)\,.$$  Furthermore,
$\hat V_{\lam\v}(\Lam_0):=\hat V_z(\Lam_0)$ for $z=t_{-\nu(\lam\v)}$ is a
$\g$-module for any dominant integral coweight $\lam\v\in
\oplus_{i=1}^r\NN\om_i\v$.

The combinatorial counterpart of this construction is:\
$$
\Bhat_{\sig y}(\pi):=\sig\Bhat_y(\pi)
$$
for $\sig\in\Sig$ and $y\in\What$.  All of our statements
regarding $\hat V_z(\Lam)$ and $\Bhat_z$ for $z\sh\in\What$ remain
valid for $z\sh\in\Wtil$.

\subsection{Proof of Theorems 1 and 3}

Let $\lam\v$ be a dominant integral coweight (not necessarily in
the root lattice) which can be written:
$$
\lam\v=\lam_1\v+\cdots+\lam_m\v\,,
$$
where $\lam_j\v\in\{\om_1,\ldots,\om_r\}$ are minuscule
fundamental coweights (not necessarily distinct) with
corresponding weights $\lam_j$ and dual weights
$\lam_j^*=-w_0(\lam_j)$.  Then we claim (extending Theorem 1 to
the coweight lattice) that there is an isomorphism of
$\g$-modules:
$$
\hat V_{\lam\v}(\Lam_0)\cong V(\lam_1^*)\otimes\cdots\otimes
V(\lam_m^*)\,.
$$
This follows immediately from Littelmann's Character and
Restriction Theorems combined with the following extension of
Theorem 3.

We will choose a certain $\ghat$-dominant path $\pi_m$ of weight
$\Lam_0$ and show that:
$$
\B_{\lam\v}(\pi_m)= \Lam_0\st
\B(\lam_1^*)\st\cdots\st\B(\lam_m^*)\,.
$$
Let $\sig_{j}\in\Sig$ correspond to $\lam_j\v$ for $j=1,\ldots,m$.
We define $\pi_m$ inductively as the last of a sequence of paths
$\pi_0,\pi_1,\ldots,\pi_m$:
$$
\pi_0:=\Lam_0,\qquad
\pi_{\!j}:=\sig_{\!j}^{\!-1}(\pi_{j-1}\st\lam_j^*)\,.
$$
We may picture $\pi_m$ as jumping up to level $\Lam_0$, winding
horizontally around the fundamental alcove
$A=(\hs_\RR\sh+\Lam_0)\cap\hat C$, and ending at $\Lam_0$. Indeed,
note that $\wt(\pi_0)=\Lam_0$.  For $j>0$, write $\om:=\lam_j$ and
$\Lam:=\om+\Lam_0$, so that $\sig_j(\Lam)=\Lam_0$ and
$\sig_j(\om)=-\om^*$.  Then by induction:
$$\begin{array}{rcl}
\wt(\pi_{j})&=&\sig_{j}^{-1}(\wt(\pi_{j-1})+\om^*)\\[.3em]
&=&\sig_{j}^{-1}\Lam_0+\sig_{j}^{-1}\om^*\\[.3em]
&=&\Lam-\om=\Lam_0\,,
\end{array}$$
so that each $\pi_j$ has weight $\Lam_0$.  Furthermore, since
$\Lam_0\sh+\lam_j^*\in\hat C$ and $\sig_j$ is a
automorphism of $\hat C$, it is clear that each $\pi_j$ is indeed
a $\ghat$-dominant path.

We will now prove Theorem 3 by showing that the Demazure operator
$\B_{\lam\v}$ ``unwinds'' $\pi_m$ starting from its endpoint. To
compute:
$$
\B_{\lam\v}(\pi_m)= \B_{\lam_1\v}\!\cdots\,
\B_{\lam_m\v}(\pi_m)\,,
$$
it suffices to prove:
\begin{lem} For $j=m,\,m\sh-1,\ldots,1$, we have:
$$
\B_{\lam_j\v}\left(\pi_{j}\st\B(\lam_{j+1}^*)
\st\cdots\st\B(\lam_m^*)\right)\hspace{.5in}
$$\\[-3em]$$\hspace{.5in}
=\pi_{j-1}\st\B(\lam_j^*)\st\B(\lam_{j+1}^*)\st\cdots\st\B(\lam_m^*)\,.
$$
\end{lem}
{\it Proof.}\ \ For $j=m$, we compute:
$$\begin{array}{rcl}
\B_{\lam_m\v}(\pi_m)&=&
\Bhat_{w_mw_0\sig_m}\left(\sig_m^{-1}(\pi_{m-1}\st\lam_m^*)\right)\\[.3em]
&=& \Bhat_{w_mw_0}\left(\pi_{m-1}\st\lam_m^*\right)\\[.3em]
&\stackrel{\mathrm{(I)}}=& \Bhat_{w_0w_m^*}\left(\pi_{m-1}\st\lam_m^*\right)\\[.3em]
&\stackrel{\mathrm{(II)}}=&
\pi_{m-1}\st\Bhat_{w_mw_0}(\lam_m^*)\\[.3em]
&\stackrel{\mathrm{(III)}}=& \pi_{m-1}\st\Bhat_{w_0}(\lam_m^*)\\[.3em]
&=& \pi_{m-1}\st\B(\lam_m^*)
\end{array}$$
The equalities are justified as follows: 
(I) Here
$w_m^*:=w_0w_mw_0$, the longest element in
$\Stab_W(\lam_m^*)=\Stab_W(-w_0\,\lam_m)$. 
(II)\ We say a path $\pi$ is {\it $i$-neutral} if
$\langle\pi(t),\al_i\v\rangle\sh\geq 0$ and 
$\langle\wt\pi,\al_i\v\rangle\sh=0$.  It is clear from the definition of $f_i$ that $f_i(\pi\st\pi')=\pi\st f_i(\pi')$ for any path $\pi'$, and the two sides are both defined or both undefined.  Now note that the $\pi_j$ are $i$-neutral for $i\sh=1,\ldots,r$.
(III)\ Follows from
$\Bhat_{w_0w_m^*}(\lam^*)=\Bhat_{w_0w_m^*}\Bhat_{w_m^*}(\lam^*)=
\Bhat_{w_0}(\lam^*)$.
\\[1em]
For $j<m$, letting $\B_{j+1}:=\B(\lam_{j+1}^*)
\st\cdots\st\B(\lam_m^*)$, we have:
$$\begin{array}{rcl}
\B_{\lam_j\v}(\pi_j\st\B_{j+1})&=&
\Bhat_{w_0w_j^*\sig_j}\left(\sig_j^{-1}(\pi_{j-1}\st\lam_j^*)\st\B_{j+1}\right)\\[.5em]
&\stackrel{\mathrm{(I)}}=&
\Bhat_{w_0w_j^*}\left(\pi_{j-1}\st\lam_j^*\st\sig_j\B_{j+1}\right)\\[.5em]
&\stackrel{\mathrm{(IV)}}=&
\Bhat_{w_0w_j^*}\left(\pi_{j-1}\st\lam_j^*\st\B_{j+1}\right)\\[.5em]
&\stackrel{\mathrm{(II)}}=&
\pi_{j-1}\st\Bhat_{w_0w_j^*}\left(\lam_j^*\st\B_{j+1}\right)\\[.5em]
&\stackrel{\mathrm{(V)}}=&
\pi_{j-1}\st\B\left(\lam_j^*\st\B_{j+1}\right)\\[.5em]
&\stackrel{\mathrm{(VI)}}=& \pi_{j-1}\st\B(\lam_j^*)\st\B_{j+1}
\end{array}$$
Here (I) and (II) are as above.
(IV) Follows from Proposition 7. (V) Follows from
\cite[Prop.~12]{LLM}, which implies that $\lam_j^*\st\B_{j+1}$ is
isomorphic to a union of Demazure crystals $\B_y(\mu)$ with $y\geq
w_m^*$. (VI) Both sides are stable under the $f_i, e_i$ for
$i=1,\ldots,r$, and they contain the same $\g$-dominant paths, hence
they are identical.

This concludes the proof of the Lemma, and hence of Theorem 3.
\QED

\section{Semi-infinite Crystals}

\subsection{Proof of Theorem 2}

Fix a minuscule coweight $\om\v$ with corresponding weight $\om$,
dual weight $\om^*$, and automorphism $\sig\in\Sig$ of order $N$.
The $N$-fold concatenation $\B(\om^*)\st\cdots\st\B(\om^*)$
contains the path $\pi_N:=
(\sig^{-N}(\om^*)\st\cdots\st\sig^{-2}(\om^*)\st\sig^{-1}(\om^*))$,
which is $\g$-dominant with weight $0$.  Thus $V(\om^*)^{\otimes
N}$ possesses a corresponding invariant vector $v_N$.

Since $t_{-N\om\v}\in\What$, the twisted Demazure module
$V_{mN\om\v}(\Lam_0)$ is an ordinary Demazure submodule of $\hat
V(\Lam_0)$.  In fact, the basic $\ghat$-module, considered as a
$\g$-module, is a direct limit of these Demazure modules:
$$
\hat V(\Lam_0)=\lim_{\longrightarrow}\ V_{N\om\v} \hookrightarrow V_{2N\om\v}
\hookrightarrow V_{3N\om\v} \hookrightarrow\cdots.
$$
By Theorem 1, the $\g$-module on the right hand side is isomorphic
to:
$$
\lim_{\longrightarrow}\ V(\om^*)^{\otimes N} \stackrel{\phi_1}\hookrightarrow
V(\om^*)^{\otimes 2N} \stackrel{\phi_2}\hookrightarrow
V(\om^*)^{\otimes 3N} \stackrel{\phi_3}\hookrightarrow\cdots
$$
for some injective $\g$-module morphisms
$\phi_1,\phi_2,\phi_3,\ldots$.

Now, if $\phi,\psi:V_1\hookrightarrow V_2$ are any two injective
$\g$-module morphisms, complete reducibility implies that there
exist automorphisms $\xi_1,\xi_2$ forming the commutative diagram:
$$\begin{array}{c@{\ }c@{\ }c}
V_1&\stackrel{\phi}\hookrightarrow&V_2\\[.2em]
\mbox{\scriptsize $\xi_1$}\!\downarrow\!\wr\ &
&\mbox{\scriptsize $\xi_2$}\!\downarrow\!\wr\ \\
V_1&\stackrel{\psi}\hookrightarrow&V_2
\end{array}$$
Thus, the inclusions $\psi_m: V(\om^*)^{\otimes Nm}\hookrightarrow
V(\om^*)^{\otimes N(m+1)}$,\ \ $v\mapsto v_N\otimes v$ appearing
in Theorem 2 define the same $\g$-module as the above maps
$\phi_m$, and we conclude that $V(\Lam_0)$ is isomorphic as a
$\g$-module to the infinite tensor product as claimed.

\subsection{The Skein model}

In this section we prove Theorem 4 by a crystal analog of the above argument.  
We introduce a notation for a path $\pi$ which emphasizes the
vector steps going toward the endpoint $\Lam=\wt(\pi)$ rather than
away from the starting point $0$.  Define
$$
\pi=(\st\, v_k\st\cdots\st v_1\ot\Lam):=(v'\!\st v_k\st\cdots\st
v_1)\,,
$$
the path with endpoint $\Lam$, last step $v_1$, etc, and first
step $v':=\Lam\sh-(v_k\sh+\cdots\sh+ v_1)$, a makeweight to
assure that the steps add up to $\Lam$.  

A {\it skein} is an infinite list:
$$
\pi=(\cdots\st v_2\st v_1\ot\Lam)\,,
$$
where $\Lam\in\oplus_{i=0}^r\ZZ\Lam_i$ and $v_j\in\hs_\RR$ (not
$\hhat^*_\RR$), subject to conditions (i) and (ii) below.  For
$i=0,\ldots,r$ and $k>0$, define:
$$
h_i[k]:=\langle\,\Lam\sh-(v_1\sh+\cdots\sh+v_k),\,\al_i\v\,\rangle\,.
$$
We require:\\
(i)\ For each $i$ and all $k\gg0$, we have $h_i[k]\geq 0$.\\
(ii)\ For each $i$, there are infinitely many $k$ such that
$h_i[k]=0$.
\\[.3em]
We think of the skein $\pi$ as a ``projective limit'' of the paths
$$
\pi[k]:=(\st\, v_k\st\cdots\st v_1\ot\Lam)\quad\text{as}\quad
k\to\infty\,.
$$
The conditions on $\pi$ assure that only a finite number of steps
of $\pi$ lie outside the fundamental chamber $\hat C$, and that
$\pi$ touches each wall of $\hat C$ infinitely many times.  Note
that $\pi$ stays always at the level $\ell=\langle \Lam,K\rangle$.

\begin{lem}
For a skein $\pi$ and $i=0,\ldots,r$, one of the following is
true:\\[.3em]
(i) $f_i(\pi[k])$ is undefined for all $k\sh\gg 0$; \\[.1em]
(ii) there is a unique skein $\pi'$ such that
$\pi'[k]=f_i(\pi[k])$ for all $k\sh\gg 0$. \\[.3em]
In the second case, we define $f_i\pi:=\pi'$.
\end{lem}
{\it Proof.}  Recall that a path $\pi$ is $i$-neutral if $h^\pi_i(t)\geq0$ for all $t$ and $h^\pi_i(1)=0$.  For a fixed $i$, divide
$\pi$ into a concatenation:
$\pi=(\cdots\st\pi_2\st\pi_1\st\pi_0\ot\Lam)$, where each $\pi_j$
is an $i$-neutral finite path except for $\pi_0$, which is an
arbitrary finite path. Now it is clear that if $f_i(\pi_0)$ is
undefined, then (i) holds. Otherwise (ii) holds and
\\[-2em]

$$
\qquad\qquad f_i\pi=(\cdots\st \pi_2\st\pi_1\st f_i(\pi_0)\,\ot\,
\Lam\sh-\al_i\,)\,.\qquad\QED
$$
\\[-1.5em]

We can immediately carry over the definitions of the path
model to skeins, including that of (Demazure) path crystals. For
example, we say that $\pi$ is an integral dominant skein if
$\pi[k]$ is integral dominant for $k\sh\gg0$, and hence for all $k$. There exist
integral dominant skeins of level $\ell=1$ only when $\g$ has a
minuscule coweight. 
We cannot concatenate two skeins, but we can concatenate a skein $\pi_1$ and a path $\pi_0$: that is,
$
\pi_1\st\pi_0:=(\,\pi_1\sh\st\pi_0\ \ot\
\wt(\pi_1)\sh+\wt(\pi_0)\,)\,.
$

\begin{prop}
For an integral dominant skein $\pi$ of weight $\Lam$,
the crystal $\Bhat(\pi)$ is a model for $\hat V(\Lam)$, and 
for Demazure modules.
\end{prop}
{\it Proof.}\ Given an integral dominant skein $\pi$ and a Weyl group element $z\in\Wtil$, we can divide $\pi=\pi_1\st\pi_0$ in such a way that the Demazure operator $\Bhat_z$ acts on $\pi$ by reflecting intervals in $\pi_0$ rather than $\pi_1$.  This gives an isomorphism between the Demazure crystals generated by the path $\wt(\pi_1)\st\pi_0$ and by the skein $\pi$:
$$\begin{array}{rcl}
\Bhat_z(\wt(\pi_1)\st\pi_0)&\stackrel\sim\to&\Bhat_z(\pi_1\st\pi_0)=\Bhat_z(\pi)\\[.2em]
\wt(\pi_1)\st\pi'&\mapsto&\pi_1\st\pi'\qquad\qquad\qquad\,.
\end{array}$$
This proves the assertion about Demazure modules.

Now, given an infinite chain of Weyl group elements $z_1\sh<z_2\sh<\cdots$, we have the morphisms of $\ghat$-crystals:
$$\begin{array}{ccccc}
\Bhat_{z_1}(\Lam)&\stackrel\sim\leftarrow&
\Bhat_{z_1}(\wt(\pi_1)\st\pi_0)&\stackrel\sim\to&\Bhat_{z_1}(\pi)\\
\cap&&&&\cap\\
\Bhat_{z_1}(\Lam)&\stackrel\sim\leftarrow&
\Bhat_{z_2}(\wt(\pi_1')\st\pi_0')&\stackrel\sim\to &\Bhat_{z_2}(\pi)\\
\cap&&&&\cap\\[-.2em]
\vdots&&&&\vdots\\
\Bhat(\Lam)&&&&\Bhat(\pi)
\end{array}$$
Since the $\ghat$ crystals at the bottom are the unions of  their Demazure crystals, they are isomorphic:  $\Bhat(\Lam)\cong\Bhat(\pi)$.\quad\QED

Now consider the situation of Theorem 4. Define the skein $\pi_\infty:=(\cdots\st\pi_N\st\pi_N\ot\Lam_0)$, and consider the commutative diagram:
$$\begin{array}{ccc}
\Bhat_{N\om\v}(\Lam_0\st\pi_N)&\stackrel\sim\to&\Bhat_{\om\v}(\pi_\infty)\\
\downarrow&&\cap\\
\Bhat_{2N\om\v}(\Lam_0\st\pi_N\st\pi_N)&\stackrel\sim\to &\Bhat_{2N\om\v}(\pi_\infty)\\
\downarrow&&\cap\\[-.2em]
\vdots&&\vdots\\
&&\Bhat(\pi_\infty)
\end{array}$$
where the vertical maps on the left are those of the direct limit.
Theorem 4 now follows from the above Proposition.

\pagebreak

\end{document}